\documentclass[11pt]{article}
\usepackage{fullpage,amssymb,latexsym}
\newenvironment{@abssec}[1]{%
     \if@twocolumn
       \section*{#1}%
     \else
       \vspace{.05in}\footnotesize
       \parindent .2in
         {\bfseries #1. }\ignorespaces
     \fi}
     {\if@twocolumn\else\par\vspace{.1in}\fi}
\newenvironment{keywords}{\begin{@abssec}{Key words}}{\end{@abssec}}
\newenvironment{AMS}{\begin{@abssec}{AMS subject classification}}{\end{@abssec}}
\newtheorem{prop}{Proposition}
\newtheorem{thm}[prop]{Theorem}
\newtheorem{lem}[prop]{Lemma}
\newtheorem{corol}[prop]{Corollary}
\newtheorem{res}[prop]{Result}

\newtheorem{ex}[prop]{Example}

\def\diag{\mathop{\rm diag}\nolimits}

\def\eqbd{\mathop{{:}{=}}}
\def\bdeq{\mathop{{=}{:}}}

\def\openC{{\rm C\kern-.48em\vrule width.06em height.6em depth-.02em 
                 \kern.48em}}
\def\openR{{{\rm I}\kern-.16em {\rm R}}}
\def\openZ{{{\rm Z}\kern-.28em{\rm Z}}}
\def\openT{{{\rm T}\kern-.42em {\rm T}}}
\def\openH{{{\rm I}\kern-.16em {\rm H}}}
\def\openK{{{\rm I}\kern-.16em {\rm K}}}
\def\openL{{{\rm I}\kern-.16em {\rm L}}}
\def\openM{{{\rm I}\kern-.16em {\rm M}}}
\def\openN{{{\rm I}\kern-.16em {\rm N}}}
\def\openP{{{\rm I}\kern-.16em {\rm P}}}

\def\eqbd{\mathop{:}{=}}
\def\ii{{\rm i}}

\let\C\openC
\let\N\openN

\let\R\openR
\let\Z\openZ
\def\nt{\noindent}
\def\prf{\nt {\bf{Proof.} {\hskip 0.07cm} }}
\def\rem{\nt {\bf{Remark.} {\hskip 0.07cm} }}
\def\eop{\hfill
        {\ \vbox{\hrule\hbox{\vrule height1.3ex\hskip0.8ex\vrule}\hrule}}
        \vskip 0.3cm \par}

\def\belowrightarrow#1{{{{}\over\ #1\ }\kern-1.1em\to}}
\def\l2{{L_2}}

\def\Mu{{{\cal M}}}

\def\tr{{\rm tr\,}}

\def\F#1{{\cal F}_{#1}}
\def\gin#1{{\lfloor {#1} \rfloor}}

\begin{document}

\title{Functions preserving nonnegativity of matrices}
\author{Gautam Bharali \\
Department of Mathematics \\
Indian Institute of Science \\
Bangalore 560012 India \and
Olga Holtz \\
Department of Mathematics \\
University of California \\
Berkeley, CA 94720 USA}
\date{November 12, 2005, revised April 5, 2007}
\maketitle

\begin{keywords} Nonnegative inverse eigenvalue problem, circulant matrices, 
(block) upper-triangular matrices, symmetric matrices, positive definite 
matrices, entire functions, divided differences. 
\end{keywords}

\begin{AMS}
15A29, 15A48, 15A42
\end{AMS}

\begin{abstract} The main goal of this work is to determine which
entire functions preserve nonnegativity of matrices of a fixed
order $n$ --- i.e., to characterize entire functions $f$ with the property 
that $f(A)$ is entrywise nonnegative for every entrywise nonnegative matrix
$A$ of size $n\times n$. Towards this goal, we present a complete characterization 
of functions preserving nonnegativity of (block) upper-triangular matrices and
those preserving nonnegativity of circulant matrices. We also derive
necessary conditions and sufficient conditions for entire functions that 
preserve nonnegativity of symmetric matrices. We also show that some of
these latter conditions characterize the even or odd functions that preserve
nonnegativity of symmetric matrices.
\end{abstract}

\section{Motivation}

The purpose of this paper is to investigate which entire functions preserve 
nonnegativity of matrices of a fixed order. More specifically, we consider 
several classes of structured matrices whose structure is preserved by entire
functions and characterize those entire functions $f$ with the property that
$f(A)$ is entrywise nonnegative for each entrywise nonnegative matrix $A$ of
size $n\times n$. The characterizations that we obtain might be of independent
interest in matrix theory and other areas of mathematics. One of our own motivations 
behind our investigation is its relevance to the inverse eigenvalue problem for
nonnegative matrices.

The long-standing inverse eigenvalue problem for nonnegative matrices is the
problem of determining, given  an $n$-tuple (multiset) $\Lambda$ of complex 
numbers, whether there exists an entrywise nonnegative matrix $A$ whose
spectrum $\sigma(A)$ is $\Lambda$. The literature on the subject is
vast and we make no attempt to review it. The interested reader is
referred to books~\cite{Minc} and~\cite{BermanPlemmons}, expository 
papers~\cite{Chu}, \cite{Eglestonetal}, \cite{Laffey1}, \cite{Laffey2}  and 
references therein, as well as to some recent papers~\cite{McDonaldNeumann}, 
\cite{ChuXu}, \cite{DuarteJohnson}, \cite{Smigoc2004}, \cite{Smigoc}, 
\cite{SotoMoro}, \cite{LaffeySmigoc}, \cite{Orsi}.

The necessary conditions for a given $n$-tuple to be realizable as the
spectrum of a nonnegative matrix known so far for arbitrary values of $n$ 
can be divided into three groups: conditions for nonnegativity of moments,
Johnson-Loewy-London inequalities, and Newton's inequalities. 

Given an $n$-tuple  $\Lambda$, its {\em moments\/} are defined as follows:
$$ 
s_m(\Lambda)\eqbd \sum_{\lambda \in \Lambda} \lambda^m ,\qquad
m \in \N. 
$$
If $\Lambda=\sigma(A)$ for some nonnegative matrix $A$, then 
$s_m(\Lambda)$ is nothing but the trace $\tr (A^m)$, and therefore must be 
nonnegative.
Another basic condition follows from the Perron-Frobenius theory~\cite{Perron},
\cite{Frobenius}:
the largest absolute value $\max_{\lambda\in \Lambda} |\lambda|$ must
be the Perron eigenvalue of a realizing matrix $A$ and therefore must itself
be in $\Lambda$. Finally, the multiset $\Lambda$ must be closed under complex
conjugation, being the spectrum of a real matrix $A$. Interestingly, the
last two conditions are in fact not independent conditions, but follow 
from the nonnegativity of moments, as was shown by Friedland in~\cite{Friedland}.
Thus, there turns out to be just one set of basic conditions
$$  
s_m(\Lambda)\geq 0\qquad {\rm for}\;\; m\in \N. 
$$
 
The next set of necessary conditions was discovered independently by 
Loewy and London in~\cite{LoewyLondon} and by Johnson in~\cite{Johnson}. 
These conditions relate moments among themselves as follows:
$$ 
s_k^m(\Lambda) \leq n^{m-1} s_{km}(\Lambda), \quad k, m\in \N. 
$$

Newton's inequalities were conjectured in~\cite{HoltzSchneider} and proved for 
$M$-matrices in~\cite{Holtz1}. An $M$-{\em{matrix}\/} is a matrix of the form 
$rI-A$, where $A$ is a nonnegative matrix, $r\geq \varrho(A)$, and where $\varrho(A)$
is the {\em spectral radius\/} of $A$:  
$$
\varrho(A)\eqbd \max_{\lambda\in \sigma(A)}
|\lambda|.
$$ If $M$ is an $M$-matrix of order $n$, then the normalized 
coefficients $c_j(M)$ of its characteristic polynomial defined by 
$$ 
\det (\lambda I - M) \bdeq \sum_{j=0}^n (-1)^j {n\choose j}  c_j(M) 
\lambda^{n-j} 
$$  
must satisfy Newton's inequalities
$$  
c_j^2(M) \geq c_{j-1}(M) c_{j+1}(M), \quad j=1, \ldots, n-1.
$$
Since the coefficients $c_j(M)$ are determined entirely by the spectrum 
of $M$, and the latter is obtained from the spectrum of a nonnegative
matrix $A$ by an appropriate shift, Newton's inequalities form yet 
another set of conditions necessary for an $n$-tuple to be realizable
as the spectrum of a nonnegative matrix. The above three sets of 
conditions --- i.e., nonnegativity of moments, Johnson-Loewy-London inequalities
and Newton's inequalities are all independent of each other but are not
sufficient for realizability of a given $n$-tuple (see~\cite{Holtz1}).

Quite a few sufficient conditions are also known (see, e.g., \cite{Suleimanova},
\cite{Laffey2}, \cite{Friedland}, \cite{Chu}) as well 
as certain techniques for perturbing or combining realizable $n$-tuples into 
new realizable $n$- or $m$-tuples (where $m\geq n$) (see, e.g.,~\cite{Soules},
\cite{Soto}, \cite{Smigoc}). Also, necessary 
and sufficient conditions on an $n$-tuple to serve as the nonzero part of the 
spectrum of some nonnegative matrix are due to Boyle and 
Handelman~\cite{BoyleHandelman}.

Finally, it follows from the Tarski-Seidenberg theorem~\cite{Tarski, Seidenberg} 
that all realizable $n$-tuples form a {\em semialgebraic set\/} 
(see also~\cite{Jacobson}), i.e.,  for any given  $n$, there exist only finitely 
many polynomial inequalities that are necessary and sufficient for an $n$-tuple 
$\Lambda$ to  be realizable as the spectrum of some nonnegative matrix $A$
(this observation was communicated to us by S.~Friedland):

Indeed, each realizable $n$-tuple $\Lambda=(\lambda_1, \ldots, \lambda_n)$
is characterized by the condition 
$$ 
\exists \;  A\geq 0 \;\; : \;\;  \det (\lambda I -A)=\prod_{j=1}^n (\lambda-\lambda_j).
$$
The last condition is equivalent to each elementary symmetric function
$\sigma_j(\Lambda)$ being equal to the $j$th coefficient of the characteristic 
polynomial of $A$ multiplied by $(-1)^j$ --- i.e., to the sum of all principal 
minors of $A$ of order $j$, for $j=1, \ldots, n$. Since the set of all nonnegative
matrices is a semialgebraic set in $n^2$ entries of the matrix and since each
sum of all principal minors of $A$ of order $j$ is a polynomial in the 
entries of $A$, the lists of coefficients of characteristic polynomials
of nonnegative matrices form a semialgebraic set, and hence the $n$-tuples
whose elementary symmetric functions match one of those lists also form 
a semialgebraic set by the Tarski-Seidenberg theorem. 

However, despite so many insights into the subject, and despite the results
obtained so far, the nonnegative inverse eigenvalue problem remains open. In fact,
the problem remains open when specialized to several important classes of 
structured matrices --- for instance, the class of entrywise nonnegative
symmetric matrices.

Note that the three sets of conditions on an $n$-tuple $\Lambda$ that we discussed 
above --- i.e., nonnegativity of moments, the Johnson-Loewy-London inequalities
and the Newton inequalities --- are necessary conditions for the realizability of
$\Lambda$ as the spectrum of a {\em symmetric} $n\times n$ matrix with nonnegative
entries (provided, of course, that all the entries of $\Lambda$ are now real). A
significant fraction of this paper will be devoted to an idea that has
relevance to the inverse eigenvalue problem for nonnegative symmetric matrices.
It is an idea that was first expressed by Loewy and London in~\cite{LoewyLondon}. 
When adapted to symmetric matrices, it  may be stated as
follows: Suppose a primary matrix function $f$ is known to map nonnegative
symmetric matrices of some fixed order $n$ into themselves. Thus $f(A)$ is 
nonnegative whenever $A$ is. Since $f(\sigma(A))=\sigma(f(A))$, both the spectrum
$\sigma(A)$ and its image under the map $f$ must then satisfy the aforementioned
conditions for realizability. This enlarges the class of necessary conditions for
the symmetric nonnegative inverse eigenvalue problem. Describing this larger class
would require knowing exactly what functions $f$ preserve nonnegativity of 
such matrices matrices (of a fixed order). Towards this end, we provide a 
characterization of all the {\em even} and {\em odd} entire functions that
preserve entrywise nonnegativity of nonnegative symmetric matrices.

Along the way, we also obtain complete characterizations of all entire
functions that preserve nonnegativity of the following classes of structured
matrices:
\begin{itemize}
\item Triangular and block-triangular matrices
\item Circulant matrices
\end{itemize}
We ought to add here that, for the above classes of structured matrices,
our results do not have a bearing on the nonnegative inverse eigenvalue
problems associated to them. In fact, the solutions of the latter problems
are quite straightforward. To be precise: an $n$-tuple $\Lambda$ is the
spectrum of an $n\times n$ triangular matrix if and only if all the
entries of $\Lambda$ are non-negative. As for circulants: the eigenvalues
of a circulant matrix $A$ are determined by its first row 
$\mathbf{a}:=[a_0 \; a_1\dots a_{n-1}]$ (see~\cite{Davis}), and in fact,
there is a constant matrix $\mathsf{W}$ (i.e., independent of $\mathbf{a}$ and
$A$) such that $\sigma(A)=\mathbf{a}\mathsf{W}$. Thus the realizable 
$n$-tuples in this case are of the form 
$\mathbf{a}\mathsf{W}, \ \mathbf{a}\in\R_+^n$. Nevertheless, we feel that
the problem of characterizing the functions that preserve nonnegativity
of the above classes of matrices can be of interest, independent of the
nonnegative inverse eigenvalue problem.

\section{Outline}

This paper is organized as follows. We make several preliminary observations 
in  Section~\ref{sec_prelim}. Before focusing attention on aspects of the
symmetric nonnegative inverse eigenvalue problem, we study the structured
matrices just discussed. In section~\ref{sec_tri}, we characterize the class
of functions preserving nonnegativity of triangular and block triangular matrices.
It turns out that these are characterized by nonnegativity conditions on their
divided differences over the nonnegative reals. Next, in Section~\ref{sec_circ},
we obtain a characterization of functions preserving nonnegativity of circulant
matrices. This characterization is quite different from that in 
Section~\ref{sec_tri} --- it involves linear combinations of function values
taken at certain non-real points of $\C$. In Section~\ref{sec_small}, we obtain
a complete characterization of the class $\F{n}$ for small values of~$n$.
 
The remainder of the paper is essentially devoted to functions that preserve 
nonnegativity of symmetric matrices. In Section~\ref{sec_spd}, we review existing
results in that direction. In particular, we discuss the restriction of 
\cite[Corollary 3.1]{MicchelliWilloughby} to entire functions, which claims
to provide a characterization of entire functions that preserve entrywise
nonnegativity of symmetric matrices of a fixed order. We point out that, while
this result is true when restricted to {\em nonnegative definite} nonnegative
symmetric matrices, the condition occurring in that result is {\em not 
sufficient} for an entire function to preserve nonnegativity of all symmetric
matrices. The {\em techniques} leading to \cite[Corollary 3.1]{MicchelliWilloughby},
however, turn out to be very useful. We use these techniques, along with some new 
ideas, to obtain necessary conditions and sufficient conditions, and 
characterizations of the {\em even} and {\em odd} entire functions
that preserve nonnegativity of symmetric matrices of a fixed order. This is the
content of Sections~\ref{sec_sym2} and~\ref{sec_sym3}. Because of a gap between 
the necessary and the sufficient conditions, which we also point out in 
Section~\ref{sec_sym2}, the results of that section do not provide a 
characterization of {\em all} functions preserving nonnegativity of symmetric 
matrices. We end the paper with a list of several open problems in 
Section~\ref{sec_next}, and suggest various approaches to their solution that 
we have not explored in this paper. 

\section{Notation} \label{sec_not}

We use standard notation $\R^{m\times n}$ for real matrices of size $m{\times}n$,
$\R_+$ for nonnegative reals and $\Z_+$ for nonnegative integers, $A\geq 0$ ($A>0$)
to denote that a matrix $A$ is entrywise nonnegative (positive), and $\sigma(A)$ to 
denote the spectrum of $A$. For $x\in\R$, we use $\gin{x}$ to denote the greatest 
integer that is less than or equal to $x$.

\section{Preliminaries} \label{sec_prelim}

The main goal of the paper is to characterize functions $f$ 
such that the matrix $f(A)$ is (entrywise) nonnegative for any nonnegative
matrix $A$ of order $n$. Since the primary matrix function $f(A)$ is
defined in accordance with values of $f$ and its derivatives on the
spectrum of $A$ (see, e.g.,~\cite[Sections 6.1, 6.2]{HornJohnson}), we want to avoid 
functions that are not differentiable at some points in $\C$. Therefore 
we restrict ourselves to functions that are  analytic everywhere in $\C$, 
i.e., to entire functions. Thus we consider the class
$$ 
\F{n}\eqbd 
\{ \; f \; {\rm entire} \; : \; A\in \R^{n\times n},\;  A\geq 0\; \Longrightarrow \; f(A)\geq 0\}.
$$

Note right away that the classes $\F{n}$ are ordered by inclusion:

\begin{lem} \label{incl}
For any $n\in \N$, $\F{n} \supseteq \F{n+1}$. 
\end{lem}

\prf Let $A$ be a nonnegative matrix of order $n$ and let
$f\in \F{n+1}$. Consider the block-diagonal matrix $B\eqbd \diag(A,0)$ 
obtained by adding an extra zero row and column to $A$. Since $f(B)=\diag(f(A),0)$,
the matrix $f(A)$ must be nonnegative. Thus $f\in \F{n}$. \eop 

Recall that any entire function can be expanded into its Taylor series around any 
point in $\C$, and that the resulting series converges everywhere (see, 
e.g.,~\cite{Conway}). We will mostly focus on Taylor series of functions
in $\F{n}$ centered at the origin.  We start with some simple observations regarding 
a few initial Taylor coefficients of such a function.

\begin{prop} \label{coefs} 
Let $f(z)=\sum_{j=0}^\infty a_j z^j$ be a function in $\F{n}$.
Then, $a_j\geq 0$ for $j=0, \ldots, n-1$.
\end{prop}

\prf For $n=1$, the statement follows from evaluating $f$ at $0$. If $n>1$ and 
$f\in \F{n}$, then evaluate the function $f$ at the matrix 
$$ 
A \eqbd \left[ \begin{array}{cccccc} 0 & 1 & 0 & \cdots & 0 & 0 \\
0 & 0 & 1 & \cdots & 0 & 0 \\ 0 & 0 & 0 & \cdots & 0 & 0 \\ 
\vdots & \vdots & \vdots & \ddots & \vdots & \vdots \\
0 & 0 & 0 & \cdots & 0 & 1 \\
0 & 0 & 0 & \cdots & 0 & 0  \end{array}    \right].
$$
Since
$$ 
f(A) = \left[ \begin{array}{cccccc} a_0 & a_1 & a_2 & \cdots & a_{n-2} & a_{n-1} \\
0 & a_0 & a_1 & \cdots & a_{n-3} & a_{n-2} \\ 0 & 0 & a_0 & \cdots & a_{n-4} & a_{n-3} \\ 
\vdots & \vdots & \vdots & \ddots & \vdots & \vdots \\
0 & 0 & 0 & \cdots & a_0 & a_1 \\
0 & 0 & 0 & \cdots & 0 & a_0  \end{array}    \right],
$$
the entries $a_0$, $\ldots$, $a_{n-1}$ of $f(A)$ must be nonnegative. 
 This finishes the proof. \eop

\begin{corol} A function $f$ is in $\F{n}$ for all $n\in \N$ if and only if
it has the form $f(z)=\sum_{j=0}^\infty a_j z^j$ with $a_j\geq 0$ for
all $j\in \Z_+$. 
\end{corol}

\prf One direction follows from Proposition~\ref{coefs}. The other direction 
is trivial: if all terms in the Taylor expansion of $f$ around the origin
are nonnegative, then $f(A)$ combines powers of a nonnegative matrix $A$
using nonnegative coefficients, so the resulting matrix is nonnegative. 
Here we make use of the standard fact~\cite[Theorem~6.2.8]{HornJohnson} that
the matrix power series $\sum_{j=0}^\infty a_j A^j$ converges to $f(A)$.  
\eop

\rem It must be noted that Proposition~\ref{coefs} {\em cannot} be a necessary 
condition for an entire function to belong to $\F{n}$. This is easy to see; fix
an $n\in \N$ and set
$$
F(x) = -x^n+\sum_{j=0}^{n-1}a_jx^j,
$$
where we choose $a_j\geq 0, \ j=0,\dots,n-1$. Then, there exists an $x_0>0$ such that 
$F(x)<0 \ \forall x\in(x_0,\infty)$. If we set $A=rI$, for some $r\in(x_0,\infty)$,
then $A$ is entrywise nonnegative while the diagonal entries of $F(A)$ are 
negative. Hence, although $a_j\geq 0$ for $j=0,\dots,n-1$, $F$ does not preserve
nonnegativity.
 
To conclude this section, we make two more general observations. 

\begin{lem} \label{poslem} An entire function $f$ belongs to $\F{n}$
if and only if it maps positive matrices of order $n$ into nonnegative 
matrices.
\end{lem} 

\prf This is simply due to the continuity of $f$, since the set of strictly positive 
matrices is dense in the set of all nonnegative matrices of order $n$. \eop

\begin{lem} \label{simlem} For any primary matrix function $f$, any permutation matrix $P$ 
and any diagonal matrix $D$ with positive diagonal elements, $f(A)$ is nonnegative 
if and only if $f(PD A (PD)^{-1})$ is nonnegative. \end{lem}

\prf  Note that $(PD)f(A)(PD)^{-1} = f(PD A (PD)^{-1})$ and that both matrices
$PD$ and $(PD)^{-1}$ are nonnegative.  So, $f(A)$ is nonnegative if and only 
if the matrix  $f(PD A (PD)^{-1})$ is nonnegative.  \eop

We now analyze three superclasses of our class $\F{n}$:
\begin{itemize}
\item entire functions preserving nonnegativity of upper-triangular
matrices;
\item entire functions preserving nonnegativity of circulant matrices; 
\item entire functions preserving nonnegativity of symmetric matrices.
\end{itemize}

\section{Preserving nonnegativity of (block-)triangular matrices} \label{sec_tri}

We first discuss functions preserving nonnegativity of upper- 
(or lower-)triangular matrices. The characterization that we obtain makes use
of the notion of divided differences. The {\em divided difference\/} (see, 
e.g.,~\cite{deBoor}) of a smooth function $f$ at points $x_1$, $\ldots$, $x_k$ 
(which can be thought of as an ordered sequence $x_1\leq \cdots \leq x_k$)
is usually defined  via the recurrence relation 
$$ 
f[x_1, \ldots, x_k]\eqbd 
\cases{\frac{f[x_2, \ldots, x_k]- f[x_1, \ldots, x_{k-1}]}{x_k-x_1} & $x_1\neq x_k$, \cr
	{} & ${ \ }$ \cr
	f^{(k-1)}(x_1)/(k-1)!  & $x_1=x_k$,} 				
$$
and where $f[x]\eqbd f(x)$. Divided differences play a large part in this paper. We shall,
however, make no attempt to review the results on divided differences that we shall draw
upon, especially since they are quite readily accessible. The interested reader is referred
to~\cite{deBoor}. 

\begin{thm} \label{thm_tri} An entire function $f$ preserves nonnegativity of 
upper-triangular matrices of order $n$ if and only if its divided differences
of order up to $n$ are nonnegative over $\R_+$, i.e., 
\begin{equation}
 f[x_1, \ldots, x_k]\ge 0 \quad {\rm for}\;\; x_1, \ldots, x_k \ge 0, \quad 
 k=1, \ldots, n, \label{divdif}
\end{equation}
or, equivalently, that all derivatives of $f$ of order up to $n{-}1$ 
are nonnegative on $\R_+$.
\end{thm}

\prf   {\em Sufficiency:} Let $A\bdeq (a_{ij})$ be a nonnegative 
upper-triangular matrix. Suppose a function $f$ satisfies~(\ref{divdif}).
By~\cite{Schmitt}, \cite{Stafney} (see also~\cite{Stafney_cor}), the elements 
of the matrix $f(A)$ can be written explicitly as
\begin{equation}
 f(A)_{ij}=  \cases{ f(a_{ii}) & $i=j$, \cr \sum_{i<i_1<\cdots<i_k<j} 
 a_{ii_1}\cdots a_{i_k j} f[a_{ii},a_{i_1i_1},\ldots, a_{i_k i_k},a_{jj}]
 & $i<j$, \cr 0 & $i>j$.}   \label{expand} 
\end{equation}
The divided differences appearing in the sum on the right-hand side are of 
order not exceeding $n$; hence all the summands, and therefore the sums, 
are nonnegative. 

{\em Necessity:} We proceed by induction on $n$. If $f$ preserves nonnegativity of 
upper-triangular matrices of order $n$, it does so also for matrices of order $n-1$. 
Thus, by our inductive hypothesis, (\ref{divdif}) holds up to order $n-1$. To see that
all divided differences of order $n$ are also nonnegative over nonnegative reals, 
consider the matrix $A$ whose first upper diagonal consists of ones, main diagonal of 
$n$ arbitrary nonnegative numbers $x_1, \ldots, x_n$, and all of whose other entries
are zero. Then,~(\ref{expand}) shows that $f(A)_{1n}=f[x_1,\ldots, x_n]$ and must be 
nonnegative.

Finally, since all divided differences of a fixed order $k$ at points in 
a domain $D$ are nonnegative if and only if $f^{(k-1)}(x)$ is nonnegative for 
every point $x\in D$~\cite{deBoor}, we see that condition~(\ref{divdif})
is equivalent to all derivatives of $f$ of order up to $n{-}1$ being nonnegative
on $\R_+$. This finishes the proof. 
\eop

The proofs of~(\ref{expand}) in~\cite{Stafney} and~\cite{Schmitt} are based on 
the following observation.  

\begin{res}[{\cite{Stafney}, \cite{Schmitt}}] \label{triang}
A block-triangular matrix of the form
$$ 
M=\left[ \begin{array}{cc} A  & B \\ 0 & a  \end{array} \right],\qquad a\in 
\C\setminus \sigma(A),
$$ 
 is mapped to the matrix 
$$ f(M)=\left[ \begin{array}{cc} f(A)  & 
(A-aI)^{-1} (f(A)-f(a)I) B \\ 
0 & f(a) \end{array}   \right] 
$$
by a function $f$. \end{res}

One can prove an analogous statement in the block-triangular case:

\begin{prop} \label{blocktriang} Let $f$ be an entire function and let
$$ 
M=\left[ \begin{array}{cc} A  & B \\ 0 & C \end{array}  \right],\qquad 
\sigma(A)\cap \sigma(C)=\emptyset.
$$  
Then,
$$ f(M)=\left[ \begin{array}{cc} f(A)  & f(A)X-Xf(C)  \\ 
0 & f(C) \end{array}  \right], 
$$
where $X$ is the (unique) solution to the equation
$$ 
AX-XC=B.  
$$
\end{prop}

\prf Let $X$ be a solution of the Sylvester equation $AX-XC=B$.
Since the spectra of $A$ and $C$ are disjoint, this solution is 
unique~\cite[Section 4.4]{HornJohnson}.  Then, $M=T^{-1}\diag(A,C)T$ where
$$  
T=\left[ \begin{array}{cc} I  & X \\ 0 & I \end{array}  \right].
$$ 
Hence $f(M)=T^{-1}\diag(f(A), f(C)) T$, which proves the proposition. 
\eop

As an immediate corollary, we obtain an indirect characterization 
of functions preserving nonnegativity of block-triangular matrices
with two diagonal blocks.

\begin{corol} An entire function  $f$ preserves nonnegativity of block 
upper-triangular matrices of the form 
$$ 
\left[ \begin{array}{cc} A  & B \\ 0 & C \end{array}  \right],\qquad
A\in \R^{n_1{\times}n_1}, \;\; C \in \R^{n_2{\times}n_2},
$$
if and only if 
\begin{itemize}
\item[a)] $f\in \F{N}$, where $N\eqbd \max\{n_1,n_2 \}$; and
\item[b)] $f(A)X-Xf(C)\geq 0$ for every $A\in \R^{n_1{\times}n_1}$,
$B\in \R^{n_1{\times}n_2}$, $C \in \R^{n_2{\times}n_2}$
such that $A, B, C \geq 0$, $\sigma(A)\cap \sigma(C)=\emptyset$, and the
(unique) matrix  $X$ that satisfies the equation $AX-XC=B$.
\end{itemize}
\end{corol}
 
\prf For $f$ to preserve nonnegativity of blocks $A$ and $C$, it has
to belong to $\F{N}$ (keeping in mind Lemma~\ref{incl}).  The remainder
of our assertion follows from Proposition~\ref{blocktriang} and the fact
that the matrices with nonnegative blocks $A$, $B$, $C$, such that the spectra
of $A$ and $C$ are disjoint, are dense in the set of all block 
upper-triangular matrices.  
\eop

The above proposition, however, does not allow for an explicit formula
of the type~(\ref{divdif}) as in Theorem~\ref{thm_tri}. \vskip 0.2cm

\rem Note that the results of this section characterize functions preserving
nonnegativity of the (block) lower-triangular matrices as well. 

\section{Preserving nonnegativity of circulant matrices} \label{sec_circ}

A {\em circulant matrix\/} (see, e.g.,~\cite{Davis}) $A$ is determined by its 
first row $(a_0, \ldots, a_{n-1})$ as follows:
$$ 
\left[ \begin{array}{ccccc} a_0 & a_1 & a_2 & \cdots & a_{n-1} \\
a_{n-1} & a_0 & a_1 & \cdots & a_{n-2} \\
a_{n-2} & a_{n-1} & a_0 & \cdots & a_{n-3} \\
\vdots & \vdots & \vdots & \ddots & \vdots \\
a_1 & a_2 & a_3 & \cdots &   a_0 \end{array}  \right]. 
$$
All circulant matrices of size $n$ are polynomials in the basic circulant matrix
$$ 
\left[ \begin{array}{ccccc} 0 & 1 & 0 & \cdots & 0 \\
0 & 0 & 1 & \cdots & 0 \\
0 & 0 & 0 & \cdots & 0 \\
\vdots & \vdots & \vdots & \ddots & \vdots \\
1 & 0 & 0 & \cdots &   0 \end{array}  \right], 
$$
which implies in particular that any function $f(A)$ of a circulant matrix
is a circulant matrix as well. Moreover, the eigenvalues of a circulant
matrix are determined by its first row~(see~\cite{Davis}) by the formula
$$
\{ \; \sum_{j=0}^{n-1} \omega^{kj} a_j \; : \; k=0, \ldots, n-1 \}, \;\;\;
{\rm where} \;\; \omega\eqbd e^{2\pi \ii/n}.
$$ 
Hence the eigenvalues of $f(A)$ are   
$$ 
\{ f( \; \sum_{j=0}^{n-1} \omega^{kj} a_j ) \; : \; k=0, \ldots, n-1 \}. 
$$
Thus, the elements $(f_0, \ldots, f_{n-1})$ of the first row of $f(A)$ 
can be read off  from its spectrum:
$$ 
f_l={1\over n} \sum_{k=0}^{n-1}  \omega^{-lk}  f( \sum_{j=0}^{n-1}
\omega^{jk}  a_j ), \;\;\; l=0, \ldots, n-1.   
$$
This argument proves the following theorem.

\begin{thm} For an entire function $f$ to preserve nonnegativity of 
circulant matrices of order $n$, it is necessary and sufficient that for
$l=0,\dots,n-1$, 
\begin{equation}
 \sum_{k=0}^{n-1}  \omega^{-lk}  f( \sum_{j=0}^{n-1}
 \omega^{jk}  a_j )\geq 0\;\;\; whenever\;\; a_j\geq 
 0, \;\; j=0,\ldots, n-1, \;\;
 where\;\; \omega=e^{2\pi\ii/n}. \label{circsums} 
\end{equation}
\end{thm}

\section{Characterization of $\F{n}$ for small values of $n$} \label{sec_small}

We now focus of the function classes $\F{n}$ for small values of $n$.
Recall the inclusion $\F{n+1}\subseteq \F{n}$ from Lemma~\ref{incl}, 
which means that all conditions satisfied by the functions from $\F{n}$ 
get inherited by the functions from $\F{n+1}$. Thus we need to
find out precisely how to strengthen the conditions that determine 
$\F{n}$ to get to the next class $\F{n+1}$. 

\subsection{The case $n=1$}

A function $f$ is in $\F{1}$ if and only if $f$ maps nonnegative reals into themselves.
While this statement is in a way a characterization in itself, if $f$ is an entire function
with {\em finitely many zeros}, we can give a description of the form that $f$ takes. For
such $f$, the proposition below serves as an alternative characterization.

\begin{prop} A function $f$ having finitely many zeros is in $\F{1}$ if and only if 
it has the form
\begin{equation}
 f(z)=g(z) \prod_{\alpha, \beta} ((z+\alpha)^2+\beta^2) \prod_{\gamma}
 (z+\gamma),  \label{case1}
\end{equation}
where the $\alpha$'s and the $\beta$'s are arbitrary reals, the $\gamma$'s
are nonnegative, and $g$ is an entire function that has no zeros in $\C$
and is positive on $\R_+$.
\end{prop}

\prf First note that since $f$ takes real values over the nonnegative reals,
all its zeros occur in conjugate pairs. Moreover, while the multiplicity
of the real negative zeros is not resticted in any way, the nonnegative zeros
must occur with even multiplicities. This produces exactly the factors
recorded in~(\ref{case1}), with nonnegative zeros corresponding to $\beta=0$. 
After factoring out all the linear factors, we are left with an entire 
function --- which we call $g(z)$ --- that has no zeros, and takes only 
positive values on $\R_+$. This gives us the expression~(\ref{case1}). \eop

\rem Incidentally, all polynomials $f$ that take only positive values on $\R_+$
are characterized by a theorem due to Poincar\'e and P\'{o}lya (see, 
e.g.,~\cite[p.175]{D'Angelo}): there exists a number $N\in \Z_+$ such that the 
polynomial $ (1+z)^N f(z)$ must have positive coefficients. Since we include 
non-polynomial functions into our class $\F{1}$, and since we allow functions to
have zeros in $\R_+$, the Poincar\'e-P\'{o}lya characterization is not directly 
relevant to our setup. 

\subsection{The case $n=2$}

We just saw that functions in $\F{1}$ are characterized by one inequality, viz.
\begin{equation} 
 f(x)\geq 0 \;\;\;\; \forall \; x\geq 0. \label{1by1} 
\end{equation}
In this subsection we will see that functions in $\F{2}$ are characterized
by two inequalities, one involving a  divided difference. We recall two
preliminary observations, Lemmas~\ref{poslem} and~\ref{simlem} that were
proved in Section~\ref{sec_prelim}. Their specialization to the case $n=2$
gives the following corollary. 
 
\begin{corol} \label{sym2} An entire function $f$ belongs to $\F{2}$ if and only if
it maps positive symmetric matrices of order $2$ into nonnegative 
matrices.
\end{corol}

\prf A strictly positive $2{\times}2$ matrix $A$ can be symmetrized by 
using the transformation $DAD^{-1}$, where $D$ is a diagonal matrix with
positive diagonal elements. Thus, Lemmas~\ref{poslem} and~\ref{simlem} 
imply that $f(A)$ is nonnegative for all strictly positive, and hence 
for all nonnegative matrices $A$ of order $2$, if and only if $f(A)$
is nonnegative for all symmetric matrices. \eop 

Now we are in a position to prove a characterization theorem for 
the class $\F{2}$. 

\begin{thm} \label{thm_2by2}
An entire function $f$ is in $\F{2}$ if and only if it satisfies 
the conditions   
\begin{eqnarray}  
 f(x+y)-f(x-y) &\geq& 0  \quad\forall\; x, y\geq 0,   \label{2by2-1} \\
 (x+y-z)f(x-y)+(z-x+y)f(x+y) &\geq& 0 \quad\forall \; x\geq z \geq 0, \; 
 y\geq x-z, \label{2by2-2}
\end{eqnarray}
or, equivalently, if $f$ satisfies~(\ref{2by2-1}) and the condition
\begin{equation} 
 (x+y)f(x-y)+(y-x)f(x+y) \geq 0 \quad \forall \; y\geq x \geq 0.  \label{2by2-2'}
\end{equation}
\end{thm}

\prf
If $f\in \F{2}$, then, in particular, $f$ preserves nonnegativity of
nonnegative circulant matrices. Thus, the conditions~(\ref{circsums}) are
necessary for $f$ to belong to $\F{2}$. Observe that 
the condition~(\ref{2by2-1}) is one of the two necessary 
conditions~(\ref{circsums}) in case $n=2$ (taking $a_0=x$ and $a_1=y$). 
Therefore, we need to check that the condition~(\ref{2by2-2}) is also 
necessary and that both together are sufficient. Then we also need to check
that conditions~(\ref{2by2-1}) and~(\ref{2by2-2}) are equivalent to 
conditions~(\ref{2by2-1}) and~(\ref{2by2-2'}). 
 
By Corollary~\ref{sym2}, we can restrict ourselves to the case
when $A$ is a positive symmetric matrix, i.e., when
$$
A=\left[ \begin{array}{cc} a_{11} & b \\ b & a_{22} 
\end{array} \right],\quad a_{11}, b, a_{22}>0. 
$$
Since the value of $f$ at $A$ coincides with the value of its interpolating 
polynomial of degree~1 with nodes of interpolation chosen at the eigenvalues 
of $A$~\cite[Sections 6.1, 6.2]{HornJohnson}, we get
$$ 
f(A)=f[r_1]I+f[r_1,r_2](A-r_1 I),
$$ 
where
$$
r_j\eqbd {a_{11}+a_{22}\over 2}+(-1)^j {\sqrt{(a_{11}-a_{22})^2+4b^2} 
\over 2}, \quad
j=1,2. 
$$ 
So, the off-diagonal entries of $f(A)$ are equal to 
$$f[r_1,r_2]b,
$$
while the diagonal entries are 
$$f[r_1,r_2](a_{jj}-r_1)+f(r_1), \quad j=1, 2.
$$
Writing 
\begin{eqnarray*} 
 x & \eqbd & {a_{11}+a_{22}\over 2}, \\
 y & \eqbd & { \sqrt{(a_{11}-a_{22})^2+4b^2} \over 2}, \\
 z & \eqbd & \min(a_{11}, a_{22}), 
\end{eqnarray*}
we see that the characterization for $\F{2}$ consists precisely of 
conditions~(\ref{2by2-1}) and~(\ref{2by2-2}).

It remains to prove that~(\ref{2by2-1}) and~(\ref{2by2-2}) are
equivalent to~(\ref{2by2-1}) and~(\ref{2by2-2'}). By simply taking $z=0$
in~(\ref{2by2-2}), we see that~(\ref{2by2-2}) implies~(\ref{2by2-2'}). 
So let us now assume~(\ref{2by2-1}) and~(\ref{2by2-2'}). We begin by stating 
a simple auxiliary fact. Taking $x=0$ and $y>0$ in~(\ref{2by2-1}) and~(\ref{2by2-2'}),
we get $f(y)\pm f(-y)\geq 0 \;\; \forall y>0$. We conclude from this that $f(y)\geq 0$
whenever $y\geq 0$ --- i.e., that $f$ satisfies~(\ref{1by1}).

First consider $y$ lying in the range $x-z\leq y \leq x$. In this case, we get
\begin{eqnarray*}
 && (x+y-z) f(x-y)+(z-x+y)f(x+y) \\ 
 && \qquad \qquad = (y-(x-z))(f(x+y)-f(x-y))+2yf(x-y)\geq 0 \qquad {\rm for}
 \;\;  x\geq z \geq 0. 
\end{eqnarray*} 
The nonnegativity of the second term above is a consequence of~(\ref{1by1}),
since $x-y$ is nonnegative in this case. Now if $y\geq x$, then~(\ref{2by2-1})
and~(\ref{2by2-2'}) simply imply that
\begin{eqnarray*}
 && (x+y-z) f(x-y)+(z-x+y)f(x+y) \\ 
 && \qquad \qquad = ((x+y)f(x-y)+(y-x)f(x+y))+z(f(x+y)-f(x-y)) \geq 0 \\
 && \qquad \qquad \quad \; {\rm for}
 \;\;  y\geq x\geq z \geq 0. 
\end{eqnarray*} 
The last two inequalities show that~(\ref{2by2-1}) and~(\ref{2by2-2'}) 
imply~(\ref{2by2-1}) and~(\ref{2by2-2}). This finishes the proof.
 \eop 

%

\section{Preserving nonnegative symmetric matrices} \label{sec_sym}

We now focus on the characterization problem for the class of entire functions that
preserve nonnegativity of symmetric matrices. We begin by recalling known facts
about functions that preserve nonnegative symmetric matrices that are in addition
nonnegative definite, i.e., have only nonnegative eigenvalues.

\subsection{Preserving nonnegative definite nonnegative symmetric matrices} 
\label{sec_spd}

Interestingly, the condition necessary and sufficient for preserving nonnegative 
symmetric matrices that are nonnegative definite  turns out to be exactly the same 
as the condition for preserving upper- (or lower-)triangular nonnegative matrices.  

The characterization of functions that preserve the class of 
{\em nonnegative definite}, entrywise nonnegative symmetric matrices is due to 
Micchelli and Willoughby~\cite{MicchelliWilloughby}. We next state a version
of their result that is useful for our purposes.

\begin{res}[version {of~\cite[Corollary 3.1]{MicchelliWilloughby}}]
\label{spdres}
An entire function $f$ preserves the class of nonnegative definite, entrywise nonnegative 
symmetric matrices if and only if all the divided differences of $f$ of order up to $n$ 
are nonnegative over $\R_+$, i.e., $f$ satisfies~(\ref{divdif}) or, equivalently, all
derivatives $f^{(j)}$ of $f$ up to order $n{-}1$ are nonnegative on $\R_+$.
\end{res}
 
The proof of Result~\ref{spdres} in~\cite{MicchelliWilloughby} relies on two facts.
The first is that $f(A)$ coincides with the interpolating polynomial of $f$, with nodes
at the eigenvalues of $A$, evaluated at $A$, i.e. that
\begin{equation}
 f(A)=f[r_1] I +f[r_1,r_2](A-r_1 I)+\cdots + f[r_1,\ldots r_n] (A-r_1 I)\cdots
 (A-r_{n-1}A).  \label{Lagrange} 
\end{equation}
The second fact is the entrywise nonnegativity of all matrix products
$$  
(A-r_1 I)\cdots (A-r_j I), \qquad j=1, \ldots n-1,  
$$
which holds under the assumption that the eigenvalues $r_1, \ldots r_n$ of $A$
are ordered 
$$ 
r_1\leq r_2 \leq  \cdots \leq r_n.  
$$ 

Observe, however, that conditions~(\ref{divdif}) are not sufficient for a function
to preserve nonnegativity of {\em all\/} nonnegative symmetric matrices. 
Indeed, let $n=2$ and let 
$$
f(x)=1+x+{1\over 2} x^2 -{2\over 3} x^3 +{1\over 4} x^4. 
$$
This function satisfies the condition~(\ref{divdif}) with $n=2$, but it maps 
the matrix 
$$ 
\left[ \begin{array}{cc} 0 & M \\ M & 0 \end{array} \right] \; ,  
$$
which is {\em not nonnegative definite}, to a matrix with negative off-diagonal
entries when $M>0$ is chosen to be sufficiently large. In fact, any $M>\sqrt{3/2}$
will produce a matrix with negative entries. 

Motivated by Result~\ref{spdres}, we would therefore like to find out what conditions
are necessary and sufficient for a function to preserve nonnegativity of a nonnegative
symmetric matrices. We begin, in the next subsection, by analyzing even and odd functions.

\subsection{Even and odd functions preserving nonnegativity of symmetric  matrices}
\label{sec_sym2} 

Using the Micchelli-Willoughby result --- i.e., Result~\ref{spdres} from the previous 
section --- and an  auxiliary  result from~\cite{Holtz2}, we shall obtain a 
characterization of even and odd functions that preserve nonnegativity of symmetric 
matrices. Our proof below will require the notion of a Jacobi matrix and that of a 
symmetric anti-bidiagonal matrix. A {\em Jacobi matrix\/} is a real, nonnegative
definite, tridiagonal symmetric matrix having positive subdiagonal entries. A 
matrix $A$ is called a {\em symmetric anti-bidiagonal matrix\/} 
if it has the form
\begin{equation}
 A= \left [ \begin{array}{ccccc} 0 & 0 & \cdots & 0 & a_n \\
 0 & 0 & \cdots & a_{n-2} & a_{n-1} \\ \vdots & \vdots & \cdot & \vdots &
 \vdots  \\0 & a_{n-2} & \cdots & 0 & 0 \\ a_{n} & a_{n-1} & \cdots & 0 & 0 
 \end{array}  \right], \quad a_1, \ldots, a_n\in \R. \label{mainform}
\end{equation}
We make use of the next two results, from~\cite{MicchelliWilloughby} and from~\cite{Holtz2}.

\begin{res}[{\cite{MicchelliWilloughby}}]   
\label{aux_jacobi}
A matrix function $f$ preserves nonnegativity of symmetric nonnegative definite matrices of 
order $n$ if and only if it maps Jacobi matrices of order $n$ into nonnegative matrices
or, equivalently, if the divided differences of $f$ up to order $n$ satisfy~(\ref{divdif})
for each ordered $n$-tuple $x_1\leq x_2 \leq \cdots \leq x_n$ of eigenvalues of a Jacobi 
matrix. \end{res}

The above result is not stated in precisely these words in \cite{MicchelliWilloughby}, 
but it is easily inferred --- it lies at the heart of the proof of 
\cite[Theorem 2.2]{MicchelliWilloughby}. In addition, we shall also need
the following result:

\begin{res}[{Corollary~3, \cite{Holtz2}}] 
\label{aux_anti}
Let $\Mu$ be a positive real $n$-tuple. Then, there exists
a Jacobi matrix that realizes $\Mu$ as its spectrum and has a symmetric 
anti-bidiagonal square root of the form~(\ref{mainform}) with all $a_j$'s
positive.
\end{res}
 
We are now in a position obtain a characterization of even and odd matrix functions 
that are of interest to us.

\begin{thm}\label{even_odd} 
An even entire function $f(z)\bdeq g(z^2)$ preserves nonnegativity of symmetric matrices
of order $n$ if and only if the divided differences of $g$ up to order $n$ are
nonnegative on $\R_+$ --- i.e., if $g$ satisfies~(\ref{divdif}). An odd function 
$f(z)\bdeq zh(z^2)$ preserves nonnegativity of symmetric matrices of order $n$ if and
only if $h$ satisfies~(\ref{divdif}).  
\end{thm}

\prf Let $f$ be even. Then, $f(z)=g(z^2)$ for some entire function $g$. If a matrix $A$ is 
entrywise nonnegative symmetric, then $A^2$ is entrywise nonnegative, symmetric, and 
nonnegative definite. By Result~\ref{spdres}, if $g$ satisfies~(\ref{divdif}), then $g(A^2)$ is 
nonnegative. To prove the converse, consider an arbitrary $n$-tuple $\Mu$ of positive numbers. 
We can think of $\Mu$ as being ordered 
\begin{equation}
 \Mu=(x_1, \ldots, x_n), \qquad x_1\leq \cdots \leq x_n. \label{mu}  
\end{equation}
By Result~\ref{aux_anti}, there exists a nonnegative 
symmetric anti-bidiagonal  matrix $A$ such that $A^2$ is a Jacobi matrix with spectrum $\Mu$. 
Then, by Result~\ref{aux_jacobi}, the divided differences of $g$ must be nonnegative when 
evaluated at the first $k$ points of $\Mu$, for each $k=1, \ldots, n$. This implies, by
the standard density reasoning, that all divided differences of $g$ must be nonnegative 
over $\R_+$.

Now let $f$ be odd. Then, $f(z)=zh(z^2)$ for some entire function $h$. If all the divided
differences of $h$ up to order $n$ are nonnegative, then by the same argument as above,
$h(A^2)$ is nonnegative for each symmetric nonnegative matrix $A$, and multiplication of
$h(A^2)$ by a nonnegative matrix $A$ produces a nonnegative matrix again. To prove the 
converse, we use induction and a technique from~\cite{MicchelliWilloughby}. 
Since $f$ has to preserve nonnegativity of symmetric matrices of order $n{-}1$ as well,
we can assume the nonnegativity of the divided differences of orders $k=1, \ldots, n{-}1$. 
To prove that the $n$th divided difference is nonnegative, let $\Mu$ be an arbitrary 
positive $n$-tuple~(\ref{mu}). As above, by Result~\ref{aux_anti}, there exists a symmetric 
 anti-bidiagonal matrix $A$ such that $A^2$ is a Jacobi matrix with spectrum $\Mu$. 
By~\cite{MicchelliWilloughby}, formula~(\ref{Lagrange}) shows that the $(1,n)$ entry of 
the function $h(A^2)$ is a positive 
multiple of $f[x_1,\ldots, x_n]$, hence the $(1,1)$ entry of the product $h(A^2) A$ is
again a positive multiple of $f[x_1, \ldots, x_n]$. Thus the $n$th divided difference
has to be nonnegative as well, which finishes the proof. \eop

This theorem provides a rather natural characterization of even and odd functions that
preserve nonnegativity of symmetric matrices in terms of their divided differences.
However, the ``natural'' idea, that the even and odd parts of any entire function that
preserves nonnegativity of symmetric functions must be also nonnegativity-preserving,
turns out to be wrong. Here is an example that illustrates why that may not be the case.

\begin{ex}\label{failurex} Let 
$$ 
f(z)\eqbd \alpha + \beta z -z^3+z^5+\gamma z^6,
$$
where $\beta>1/4$, and $\alpha,\gamma>0$ are chosen to be so large
that $f(x)\geq 0$ for all $x\in \R$ and $f'(x)\geq 0$ for all $x\in \R_+$.
Then, $f$ preserves nonnegativity of symmetric matrices of order $2$,
but its odd part $f_{odd}$ does not.   
\end{ex}

\prf The function $f$ satisfies conditions~(\ref{2by2-1}) and~(\ref{2by2-2'}).
Indeed, since $f\geq 0$ on $\R$, we have
$$ 
(t+s) f(-t)+tf(t+s) \geq 0 \qquad \forall \; s,t\geq 0,
$$
which is equivalent to condition~(\ref{2by2-2'}). Now, the odd part of
$f$ is given by
$$ 
f_{odd}(z)=\beta z-z^3+z^5\bdeq zh(z^2).
$$
Since $\beta>1/4$, $h(x)>0$ for all $x\in \R$. Since $f$ is monotone
increasing on $\R_+$, we have
$$ 
f(s+t)-f(-s)\geq f(s)-f(-s)=2f_{odd}(s)\geq0 \qquad \forall\; s, t \geq 0, 
$$
which yields condition~(\ref{2by2-1}). Thus, by Theorem~\ref{thm_2by2}, $f$ preserves 
nonnegativity of symmetric matrices of order $2$. However, 
$$
h'(x)=2x-1<0 \qquad {\rm for} \;\; x<1/2.
$$
Therefore, by Theorem~\ref{even_odd}, $f_{odd}$ does not preserve 
nonnegativity of symmetric functions of order $2$. \eop

We conclude this section with a simple observation about even and odd parts
of a nonnegativity-preserving function.

\begin{prop} If an entire function $f$ preserves nonnegativity of symmetric
functions of order $n$, then its odd and even parts $f_{odd}$ and $f_{even}$
preserve nonnegativity of matrices of order $\lfloor n/2 \rfloor$.
\end{prop}

\prf For $n$ even, consider matrices of the form
$$ 
A=\left[ \begin{array}{cc} 0 & B \\ B & 0 \end{array}  \right], 
$$
and for $n$ odd, matrices of the form
$$ 
A=\diag (\left[ \begin{array}{cc} 0 & B \\ B & 0  \end{array}  \right], 0), 
$$
where $B$ is an $\lfloor n/2 \rfloor{\times}\lfloor n/2 \rfloor$
symmetric nonnegative matrix. Since 
$$ 
f(A)=\left[ \begin{array}{cc} f_{even}(B) & f_{odd}(B) \\ f_{odd}(B) & 
f_{even}(B) \end{array}  \right] \qquad {\rm for} \; n \; {\rm even}, 
$$
$$ 
f(A)=\diag(\left[ \begin{array}{cc} f_{even}(B) & f_{odd}(B) \\ f_{odd}(B) & 
f_{even}(B) \end{array}  \right], 0 ) \qquad {\rm for} \; n \; {\rm odd}, 
$$
the see that $f_{even}$ and $f_{odd}$ must preserve nonnegativity of
symmetric functions of order $\lfloor n/2 \rfloor$. \eop

\subsection{Other necessary conditions} \label{sec_sym3}

Results from~\cite{Holtz2} allow us to derive an additional set of 
necessary conditions. The motivation behind these conditions is as 
follows. We believe that the power of Results~\ref{aux_jacobi} 
and~\ref{aux_anti} --- or rather, {\em the methods} behind those
results --- have not been exhausted by Theorem~\ref{even_odd}. Our
next theorem is presented as an illustration of this viewpoint. On
comparison with Theorem~\ref{thm_2by2}, we find that the conditions
derived in our next theorem constitute a complete characterization
for the functions of interest in the $n=2$ case. To derive these 
new necessary conditions, we will need the following two results.

\begin{res}[{Theorem~1, \cite{Holtz2}}]
\label{aux_antiChar}
A real n-tuple $\Lambda$ can be realized as the spectrum of a symmetric
anti-bidiagonal matrix~(\ref{mainform}) with all $a_j$'s positive if and
only if $\Lambda=(\lambda_1,\ldots,\lambda_n)$ satisfies
$$
\lambda_1 > -\lambda_2 > \lambda_3 > \cdots > (-1)^{n-1}\lambda_n > 0.
$$
\end{res}

\begin{lem}\label{aux_pattern} Let $A$ be a symmetric anti-bidiagonal matrix
of order $n$, and  let $A^p_{ij}$ denote the $(i,j)$ entry of $A^p$. Then
\begin{enumerate}
\item[a)] The $(i,j)$ entry of $A^{2q-1}$ is zero whenever 
$2\leq i+j\leq (n-q+1)$, $q\geq 1$.
\item[b)] The $(i,j)$ entry of $A^{2q}$ is zero whenever 
$1+q\leq j-i\leq n-1$, $q\geq 1$.
\item[c)] Adopting the notation in~(\ref{mainform}) for the entries of $A$,
\begin{eqnarray}
 A^{2q-1}_{1,n-q+1} &=& a_n a_{n-1}\ldots a_{n-2q+2}, \qquad 1\leq q\leq 
 \gin{(n+1)/2}, 
 \label{prod1}\\
 A^{2q}_{1,1+q} &=&  a_n a_{n-1}\ldots a_{n-2q+1}, \qquad 1\leq q\leq \gin{n/2}.
 \label{prod2}
\end{eqnarray}
\end{enumerate}
\end{lem}

\prf We proceed by induction on $q$. Note that (a), (b) and (c) are obvious when
$q=1$. Let us now assume that (a) and (b) are true for some $q<n-3$. Note that
since $A$ is anti-bidiagonal,
\begin{equation}
 A^{2q+1}_{ij} = A^{2q}_{i,n-j+1}A_{n-j+1,j}+A^{2q}_{i,n-j+2}A_{n-j+2,j}. 
 \label{expan(a)}
\end{equation}
However, if $i+j\leq(n-(q+1)+1)$, then
$$
(n-j+2)-i\geq (n-j+1)-i\geq q+1.
$$
Applying our inductive hypothesis on (b), we conclude from the above inequalities
that the right-hand side of~(\ref{expan(a)}) reduces to zero when
$i+j\leq(n-(q+1)+1)$. Thus, (a) is established for $q+1$.

We establish (b) for $q+1$ in a similar fashion. We note that
\begin{equation}
 A^{2q+2}_{ij} = A^{2q+1}_{i,n-j+1}A_{n-j+1,j}+A^{2q+1}_{i,n-j+2}A_{n-j+2,j}. 
 \label{expan(b)}
\end{equation}
When $j-i\geq 1+(q+1)$, then
$$
i+(n-j+1)\leq i+(n-j+2)\leq n-(q+1)+1.
$$
Since we just established (a) for $q+1$, the above inequalities tell us that
the right-hand side of~(\ref{expan(b)}) reduces to zero when $j-i\geq 1+(q+1)$.
Thus, (b) too is established for $q+1$. By induction, (a) and (b) are true
for all relevant $q$.

Part (c) now follows easily by substituting $i=1$ and $j=n-q$ into
equation~(\ref{expan(a)}) to carry out the inductive step for~(\ref{prod1}),
and by substituting $i=1$ and $j=q+2$ into equation~(\ref{expan(b)}) to
carry out the inductive step for~(\ref{prod2}). 
\eop

We can now present the aforementioned necessary conditions.
 
\begin{thm}\label{thm_newnc} If an entire function $f$ preserves nonnegativity
of symmetric matrices of order $n$, $n\geq 2$, then, for each ordered 
$n$-tuple $(x_1, \ldots, x_n)$ where
\begin{equation}
 x_1 > -x_2 > x_3 > \cdots > (-1)^{n-1} x_n > 0, \label{cond}
\end{equation}
$f$ must satisfy
\begin{equation}
 f[x_1, \ldots, x_n] \geq 0,  \label{newnc1}
\end{equation}
and for each $k=1,\dots,n$, $f$ must satisfy
\begin{equation}
 f[x_1,\ldots, x_{k-1},x_{k+1}, \ldots, x_n]-( \; \sum_{j\neq k}x_j)f[x_1, \ldots, x_n]
 \geq 0. \label{newnc2}
\end{equation}
\end{thm}

\prf We choose an $n$-tuple $(x_1,\ldots,x_n)$ that satisfies~(\ref{cond}). By
Result~\ref{aux_antiChar}, there is a symmetric anti-bidiagonal matrix of the
form~(\ref{mainform}), with all $a_j$'s positive, whose spectrum is 
$(x_1,\ldots,x_n)$. Let us express $f(A)$ using the formula~(\ref{Lagrange}),
with the substitutions $r_j=x_j$, $j=1,\ldots,n$. Then, in view of 
Lemma~\ref{aux_pattern}, the $(1,\gin{n/2}+1)$ entry of $f(A)$
is $a_n a_{n-1}\ldots a_2f[x_1,\ldots,x_n]$. Since $f$ preserves
nonnegativity, and all the $a_j$'s are positive, $f[x_1,\ldots,x_n]$ has 
to be nonnegative. This establishes~(\ref{newnc1}).

To demonstrate~(\ref{newnc2}), we look at the entries of $f(A)$ that are
{\em adjacent} to the $(1,\gin{n/2}+1)$ entry that was considered above. Let us
fix a $k=1,\ldots,n$. This time, however, in using formula~(\ref{Lagrange})
to express $f(A)$, we make the following substitutions
$$
r_j = \cases{x_j & if $j<k$, \cr
	 x_{j+1} & if $k\leq j < n$, \cr
	 x_k & if $j=n$.}
$$
Our analysis splits into two cases.

{\em {\bf Case 1.} $n$ is odd:} In this case, let us look at the 
$(1,\gin{n/2}+2)$ entry of $f(A)$. By Lemma~\ref{aux_pattern}, and the fact that
$n$ is odd, the only power of $A$ that contributes to this entry is $A^{n-2}$. 
Consequently
\begin{eqnarray*}
 f(A)_{1,\gin{n/2}+2} &=&
 \{f[x_1,\ldots,x_{k-1},x_{k+1},\ldots,x_n]-( \; \sum_{j\neq k}x_j)f[x_1,\ldots,x_n] \}
 A^{n-2}_{1,\gin{n/2}+2}  \\
 &=& \ a_n a_{n-1}\ldots a_3
 \{f[x_1,\dots,x_{k-1},x_{k+1}\dots,x_n]-( \; \sum_{j\neq k}x_j)f[x_1,\ldots,x_n] \}.
\end{eqnarray*}
Since $f$ preserves nonnegativity,~(\ref{newnc2}) follows from the above equalities.

{\em {\bf Case 2.} $n$ is even:} In this case, we focus on the
$(1,\gin{n/2})$ entry of $f(A)$. We recover~(\ref{newnc2}) by arguing exactly
as above. 

In either case,~(\ref{newnc2}) is established, which concludes our proof.
\eop
 
We conclude this section by showing that a subset of the necessary conditions
derived above are in fact sufficient to characterize those entire functions
that preserve nonnegativity of $2\times 2$ symmetric matrices. Specifically,
we show that 
\begin{eqnarray*}
 f[x_1,x_2] &\geq& 0 \quad{\rm and} \\
 f(x_2)-x_2f[x_1,x_2] &\geq& 0 \quad{\rm for} \; {\rm all} \; \; x_1>-x_2>0,
\end{eqnarray*}
imply the conditions~(\ref{2by2-1}) and~(\ref{2by2-2'}). This is achieved
simply by taking some $y>x>0$, making the substitutions $x_1=y+x$ and
$x_2=x-y$, and then invoking continuity to 
obtain~(\ref{2by2-1}) and~(\ref{2by2-2'}) for all $y\geq x\geq 0$.
 
\section{Open problems and further ideas} \label{sec_next}

We conclude this paper by listing some ideas that we did not pursue, which
however may lead to further progress.

One can consider matrices that preserve nonnegativity of other classes of structured 
matrices, such as Toeplitz or Hankel. However, since these classes are not invariant
under the action of an arbitrary matrix function, their matrix functions can be quite difficult 
to analyze. Also, the eigenstructure of some structured matrices is rather involved,
which could be an additional obstacle.

Theorem~1.3 of~\cite{Stafney} gives an interesting formula for $f(A)$
when $f$ is a polynomial, which therefore must also be true for entire functions.
Precisely, if $A$ is a matrix with minimal polynomial $p_0$ and $C$ is
the companion matrix of $p_0$, then 
$$ f(A)=\sum_{j=1}^n f(C)_{j1} A^{j-1}.$$
In particular, $f(A)$ is nonnegative whenever the first column of $f(C)$
is nonnegative. It would be worthwhile to  find out what functions have this property.

Note that the set $\F{n}$ contains positive constants and is closed under 
addition, multiplication, and composition. We are not aware of any work on 
systems of entire functions (or even  polynomials) that satisfy this property. 
Perhaps one could describe a minimal set of generators (with respect to these three 
operations) that generate such a system. 

For example, in the case $n=1$, the generators are positive constants, the
function $p_1(x)=x$ plus all quadrics of the form $(x-a)^2$, $a>0$. Incidentally,
the set of polynomials with nonnegative coefficients is generated by positive
constants and $p_1(x)=x$. We do not have a characterization of generators for
$n\geq 2$.

In particular, $\F{n}$ is a semigroup with respect to any of these operations, 
so some general results on semigroups may prove to be useful in our setting.
Also note that the set of nonnegative matrices of order $n$, on which $\F{n}$ acts, 
is also a semigroup (closed under addition and multiplications), which could also 
be of potential use.

Finally, both $\F{n}$ and the set of nonnegative matrices of order $n$ are also 
cones, so the problem might also have a cone theoretic form. If we consider polynomials
instead of entire functions, we can further restrict ourselves to polynomials
of degree bounded by a fixed positive integer. Then, we will obtain a proper cone,
whose extreme directions may be of interest. The general problem then can also be looked 
upon in an appropriate similar setting.

\section*{Acknowledgments} We are grateful to Raphael Loewy, Michael Neumann
and Shmuel Friedland for helpful discussions and to anonymous referees for 
useful suggestions.

\bibliographystyle{plain}
\bibliography{matrix}

\end{document}